\documentclass[12pt]{article}
\usepackage{graphicx}
\usepackage{graphics}
\usepackage{amsthm}
\usepackage{amssymb, amsmath}
\title{Summatory relations and prime products for the Stieltjes constants, and other related results} 
\author{Mark W. Coffey\\
Department of Physics\\
Colorado School of Mines\\
Golden, CO  80401\\
(Received $\mbox{~~~~~~~~~~~~~~~~~~~~~~~~~~~~~~~2016}$)}
\date{July 31, 2016}
\begin{document}
\maketitle
\baselineskip=25 pt
\begin{abstract}

The Stieltjes constants $\gamma_k(a)$ appear in the regular part of the Laurent expansion for the Hurwitz zeta function $\zeta(s,a)$.  We present summatory results for these constants $\gamma_k(a)$ in terms of fundamental mathematical constants such as the Catalan constant, and further
relate them to products of rational functions of prime numbers.  
We provide examples of infinite series of differences of Stieltjes constants evaluating as
volumes in hyperbolic $3$-space.  We present a new series representation for the difference
of the first Stieltjes constant at rational arguments.  
We obtain expressions for $\zeta(1/2)L_{-p}(1/2)$, where for primes $p>7$, $L_{-p}(s)$ are certain
$L$-series, and remarkably tight bounds for the value $\zeta(1/2)$, $\zeta(s)=\zeta(s,1)$ being
the Riemann zeta function.




\end{abstract}
 
\medskip
\baselineskip=15pt
\centerline{\bf Key words and phrases}
\medskip 

\noindent

Dirichlet $L$ function, Hurwitz zeta function, Stieltjes constants, series representation, Riemann zeta function, Laurent expansion, digamma function, summatory relation, prime product

\vfill
\centerline{\bf 2010 AMS codes} 
11M06, 11Y60, 11M35

\baselineskip=25pt
\pagebreak
\medskip
\centerline{\bf Introduction}
\medskip

Let $\zeta(s,a)$ be the Hurwitz zeta function, and $\zeta(s,1)=\zeta(s)$ be the Riemann
zeta function \cite{edwards,ivic,riemann,titch} .  In the complex plane of $s$, each of these functions has a simple pole of residue $1$ at $s=1$.  For Re $s>1$ and Re $a>0$ we have
$$\zeta(s,a)=\sum_{n=0}^\infty {1 \over {(n+a)^s}}, \eqno(1.1)$$
and by analytic continuation $\zeta(s,a)$ extends to a meromorphic function through out
the whole complex plane. Accordingly, there is the Laurent expansion in terms of the
Stieltjes constants $\gamma_k(a)$,
\cite{berndt,coffeyjmaa,coffeystdiffs,coffey2009,stieltjes,wilton2}
$$\zeta(s,a)={1 \over {s-1}}+\sum_{n=0}^\infty {{(-1)^n} \over {n!}}\gamma_n(a)(s-1)^n,
~~~~~~s \neq 1.  \eqno(1.2)$$
We let $\Gamma(s)$ be the Gamma function, $\psi(s)=\Gamma'(s)/\Gamma(s)$ be the digamma function (e.g., \cite{nbs,andrews,grad}) with the Euler constant $\gamma=\gamma_0(1)=-\psi(1)$, and recall that $\gamma_0(a)=-\psi(a)$. In the following, $\psi^{(k)}$ are the polygamma functions.
We note that for integers $n>1$, $\zeta(n,a)$ reduces to values $\psi^{(n-1)}$.

The sequence $\{\gamma_k(a)\}_{k=0}^\infty$ exhibits complicated changes in
sign with $k$.  For instance, for both even and odd index, there are infinitely
many positive and negative values.  Furthermore, there is sign variation with the
parameter $a$.  These features, as well as the exponential growth in magnitude
in $k$, are now well captured in an asymptotic expression (\cite{knessl2}, Section 2, \cite{paris}).
In fact, though initially derived for large values of $k$, this expression is
useful for computational approximation even for small values of $k$.

Reciprocity and other relations for the Stieltjes constants, Bernoulli polynomials, and
functions $A_k(q)$ are given in \cite{coffeyrecip}.  
The works \cite{coffeygammaratl} and \cite{coffey2009} present explicit expressions for
the initial values of the Stieltjes constants at rational argument, and various series
representations, respectively.  


This paper is largely illustrative, being a step in a program to explore summatory
relations for the Stieltjes constants in terms of fundamental mathematical constants such 
as the Catalan constant
$$G \equiv \sum_{n=0}^\infty {{(-1)^n} \over {(2n+1)^2}}.  \eqno(1.3)$$  
We also relate sums of differences of Stieltjes constants to other mathematical constants 
and products of prime numbers.  
\footnote{For such prime product results, especially see Proposition 7}.
For the latter purpose, $p$ is reserved for the product index
over primes.  Two of many examples to be presented are then
$$\sum_{k=0}^\infty {1 \over {k!}}\left[\gamma_k\left({1 \over 3}\right)-\gamma_k\left({2 \over 3}\right)\right]={1 \over 3}, \eqno(1.4a)$$
and
$$G={\pi^2 \over 8}\prod_{\overset{p \equiv 3}{mod ~4}} {{p^2-1} \over {p^2+1}}={1 \over {16}}
\sum_{k=0}^\infty {1 \over {k!}}\left[\gamma_k\left({1 \over 4}\right)-\gamma_k\left({3 \over 4}\right)\right]. \eqno(1.4b)$$
While, for instance, the first equality in (1.4b) is known, it seems unlikely that the
connection with the Stieljtes constants in the second equality has been established before.


A major tool for the first sections of this paper are quadratic Dirichlet-$L$ series and
their properties.  Therefore, some further notation and definitions are introduced.  Let $D$
be a fundamental discriminant and $(D/n)$ the Kronecker-Jacobi-Legendre symbol.
This symbol, a completely multiplicative function on the positive integers, is a real primitive Dirichlet character with modulus $|D|$.  Then the Dirichlet
$L$ series associated to $(D/n)$ is defined for Re $s>1$ as
$$L_D(s)=\sum_{n=1}^\infty \left({D \over n}\right)n^{-s}. \eqno(1.5)$$
If $D=1$, then $L_1(s)=\zeta(s)$.  For all other values of $D$, $L_D(s)$ can be made into an
entire function by using the value $L_D(1) \neq 0$.  
With suitable factors of $D^{s/2}$, $\pi^{-s/2}$, and $\Gamma[(s+\epsilon)/2]$ ($\epsilon=0,1$),
$L_D(s)$ may be completed to $L^*_D(s)$ with the compact functional equation 
$L^*_D(s)=L_D^*(1-s)$.  A useful compilation of special values of $L_D(s)$ is contained in
\cite{mathar}.

Importantly for what follows, we describe how Dirichlet $L$-functions $L_{\pm k}(s)$ (e.g., \cite{ireland}, Ch. 16), are expressible as linear combinations of Hurwitz zeta functions.  
We let $\chi_k$ be a real Dirichlet character modulo $k$, where the corresponding $L$ 
function is written with subscript $\pm k$ according to $\chi_k(k-1)=\pm 1$.  We have
$$L_{\pm k}(s)=\sum_{n=1}^\infty {{\chi_k(n)} \over n^s}={1 \over k^s}\sum_{m=1}^k 
\chi_k(m) \zeta\left(s,{m \over k}\right), ~~~~\mbox{Re} ~s >1.  \eqno(1.6)$$
This equation holds for at least Re $s>1$.  If $\chi_k$ is a nonprincipal character, 
as we typically assume in the following, then convergence obtains for Re $s>0$.

These $L$ functions, extendable to the whole complex plane, satisfy the functional
equations \cite{zucker76}
$$L_{-k}(s)={1 \over \pi}(2\pi)^s k^{-s+1/2}\cos\left({{s \pi} \over 2}\right)\Gamma
(1-s)L_{-k}(1-s), \eqno(1.7)$$
and
$$L_{+k}(s)={1 \over \pi}(2\pi)^s k^{-s+1/2}\sin\left({{s \pi} \over 2}\right)\Gamma
(1-s)L_{+k}(1-s). \eqno(1.8)$$
Owing to the relation
$$\Gamma(1-s)\Gamma(s)={\pi \over {\sin \pi s}}, \eqno(1.9)$$
these functional equations may also be written in the form
$$L_{-k}(1-s)=2(2\pi)^{-s} k^{s-1/2}\sin\left({{\pi s} \over 2}\right) \Gamma(s)
L_{-k}(s), \eqno(1.10)$$
and
$$L_{+k}(1-s)=2(2\pi)^{-s} k^{s-1/2}\cos\left({{\pi s} \over 2}\right) \Gamma(s)
L_{+k}(s). \eqno(1.11)$$


%

\medskip
\centerline{\bf Statement of results}
\medskip

The first result may be placed in the context of special function theory, while the
subsequent ones Propositions 2--8 and 10 provide summatory results for differences of the Stieltjes constants.  Proposition 9 concerns representation of Sierpinski's constant $S$,
and the last Proposition 11 bounds $\zeta(s)$ in the critical strip for real values of $s$.
\newline{\bf Proposition 1} (a) 
$$L_{\pm k}(1)={1 \over k}\sum_{m=1}^k\chi_k(m)\gamma_0\left({m \over k}\right),$$
and (b) for all $D \neq 1$, the values
$$L_D(1)= \begin{cases}
{\pi \over {3\sqrt{3}}} ~\mbox{if}~ D=-3 \\
{\pi \over 4} ~\mbox{if}~ D=-4 \\
{{\pi h(D)} \over \sqrt{-D}} ~\mbox{if}~ D<-4 \\
{{2 h(D)\ln \varepsilon} \over \sqrt{D}} ~\mbox{if}~ D>1 \\
\end{cases}, $$ 
may be written in terms of $\gamma_0=-\psi$. 
Here $h(D)$ is the ideal class number of the quadratic field $\mathbb{Q}(\sqrt{D})$ and
$\varepsilon$ is the fundamental unit of the integer subring $\mathbb{Z}+((D+\sqrt{D})/2)\mathbb{Z}$
.

{\bf Proposition 2} (a) 
$$L_{\pm k}(2)={1 \over k^2}\sum_{n=0}^\infty {{(-1)^n} \over {n!}} \sum_{m=1}^k\chi_k(m)\gamma_n\left({m \over k}\right),$$
and (b) for all $D$, the values $L_D(2)$ 
may similarly be written in terms of sums of differences $\gamma_k(a)$. 
This includes the cases $L_1(2)=\zeta(2)=\pi^2/6$, $L_{-4}(2)=G$, Catalan's constant,
$$L_5(2)={{4\pi^2} \over {25\sqrt{5}}}, ~~~~L_8(2)={\pi^2 \over {8\sqrt{2}}}, ~~~~
L_{12}(2)={\pi^2 \over {6\sqrt{3}}},$$
and
$$L_{-7}(2)=I_{-7} \equiv {{24} \over {7 \sqrt{7}}}\int_{\pi/3}^{\pi/2} \ln \left|{{\tan t+\sqrt{7}} \over {\tan t-\sqrt{7}}}\right |dt. \eqno(2.1)$$

The value (2.1) and related integrals arise in hyperbolic geometry, knot theory, and quantum field
theory, and (2.1) has received considerable attention in the last several years \cite{bb07,bbnotices,bbcomp,bb98}.
References \cite{coffeyjmp,coffeyjmp2} provide alternative evaluations of the integral (2.1)
in terms of the Clausen function Cl$_2(\theta)$, 
$$\mbox{Cl}_2(\theta)\equiv -\int_0^\theta\ln\left|2 \sin {t \over 2}\right|dt
=\int_0^1 \tan^{-1}\left({{x \sin \theta} \over {1-x\cos \theta}}\right)
{{dx} \over x}$$
$$=-\sin \theta \int_0^1 {{\ln x} \over {x^2-2x\cos \theta +1}}dx
=\sum_{n=1}^\infty {{\sin(n \theta)} \over n^2}.$$
The evaluations of Proposition 2 provide connections of infinite sums of differences of Stieltjes constants with values Cl$_2(\theta)$.

The next result, among other relations, connects a sum of differences of Stieltjes constants with
a sum involving the sum of divisors function $\sigma_{-3}(n)=\sum_{d|n}d^{-3}$.  In stating the
result, we introduce the Bloch-Wigner dilogarithm
$$D(z)=\mbox{Im}[\mbox{Li}_2(z)+\ln |z|\ln(1-z)]= \mbox{Im}[\mbox{Li}_2(z)]+\mbox{arg}(1-z) \ln |z|,$$
where, as usual for $|z|\leq 1$, the dilogarithm function is given by
$\mbox{Li}_2(z)=\sum_{n=1}^\infty {z^n \over n^2}$, which analytically continues to $\mathbb{C}
\backslash (1,\infty)$. 
Then $D(z)$ is real analytic in $\mathbb{C}\backslash \{0,1\}$.
\newline{\bf Proposition 3} Let $z_1=(1+i\sqrt{23})/2$, $z_2=2+i\sqrt{23}$, $z_3=(3+i\sqrt{23})/2$, $z_4=
(5+i\sqrt{23})/2$, and $z_5=3+i\sqrt{23}$.  Then
$${3 \over 2}\zeta(4)+{{40\pi} \over 23^{3/2}}\left[{{\zeta(3)} \over 2}+\sum_{n=1}^\infty
\sigma_{-3}(n)(1+\pi n\sqrt{23})e^{-\pi n\sqrt{3}}\right]$$
$$={{4\pi^2/3} \over {23^{3/2}}}[21D(z_1)+7D(z_2)+D(z_3)-3D(z_4)+D(z_5)]$$
$$=23^{-2}\sum_{n=0}^\infty {{(-1)^n} \over {n!}}\sum_{j=1}^{22}\chi_{-23}(j)\gamma_n\left({j \over
{23}}\right)$$
$$\equiv 23^{-2}\sum_{n=0}^\infty {{(-1)^n} \over {n!}}\left[\gamma_n\left({1 \over {23}}\right)+
\gamma_n\left({2 \over {23}}\right)+\gamma_n\left({3 \over {23}}\right)+
\gamma_n\left({4 \over {23}}\right)-\gamma_n\left({5 \over {23}}\right)+
\gamma_n\left({6 \over {23}}\right) \right.$$
$$-\gamma_n\left({7 \over {23}}\right)+
\gamma_n\left({8 \over {23}}\right)+\gamma_n\left({9 \over {23}}\right)-
\gamma_n\left({{10} \over {23}}\right)+\gamma_n\left({{11} \over {23}}\right)+
\gamma_n\left({{12} \over {23}}\right)-\gamma_n\left({{13} \over {23}}\right)+\gamma_n\left({{14} \over {23}}\right)$$
$$\left. -\gamma_n\left({{15} \over {23}}\right)+\gamma_n\left({{16} \over {23}}\right)
-\gamma_n\left({{17} \over {23}}\right)+\gamma_n\left({{18} \over {23}}\right)-
\gamma_n\left({{19} \over {23}}\right)+\gamma_n\left({{20} \over {23}}\right)-\gamma_n\left({{21} \over {23}}\right)-\gamma_n\left({{22} \over {23}}\right)\right].$$

{\bf Proposition 4} (a) 
$$L_{\pm k}(3)={1 \over k^3}\sum_{n=0}^\infty {{(-1)^n} \over {n!}}2^n \sum_{m=1}^k\chi_k(m)\gamma_n \left({m \over k}\right),$$
and (b) for all $D$, the values $L_D(2)$ 
may similarly be written in terms of sums of differences $\gamma_k(a)$. 
In particular, for $D<0$, closed-form expressions for $L_D(3)$ are known.  

After the proof of Proposition 4 we illustrate the use of an integral representation for
$\gamma_k(a)$ in order to obtain such sums for rational values of $a$.

{\bf Proposition 5} 
$$4^{-4} \sum_{n=0}^\infty{{(-1)^n} \over {n!}}3^n\left[\gamma_n\left({1 \over 4}\right)-
\gamma_n\left({3 \over 4}\right)\right]={1 \over {768}}\psi'''\left({1 \over 4}\right)-{\pi^4 \over {96}}.$$

The next sort of result follows from consideration of the Euler-Kronecker constant $\gamma_D$
of the quadratic field $\mathbb{Q}(\sqrt{D})$.  It provides a new summation representation of
the difference of the first Stieltjes constant at rational arguments.
\newline{\bf Proposition 6} (a)
$$\gamma_1\left({3 \over 4}\right)-\gamma_1\left({1 \over 4}\right)={\pi^2 \over 3}+\pi \gamma
+4\pi \sum_{\ell=1}^\infty {1 \over \ell}{1 \over {(e^{2\pi \ell}-1)}},$$
and (b)
$$\gamma_1\left({2 \over 3}\right)-\gamma_1\left({1 \over 3}\right)={\pi \over \sqrt{3}}\left[
{\pi \over {2\sqrt{3}}}+\gamma -4\sum_{\ell=1}^\infty {{(-1)^\ell} \over \ell}{1 \over {[(-1)^\ell-e^{\sqrt{3}\pi \ell}]}}\right].$$

As we mention below, such differences of $\gamma_1(a)$ at rational argument are related to 
values of $\ln \Gamma(a)$ and hence $\ln \Gamma(1-a)$.  For further details \cite{coffeygammaratl}
may be consulted.

{\bf Proposition 7} (a)
$$\prod_{\overset{p \equiv 1}{mod ~3}} {{p^2+1} \over {p^2-1}}={3 \over {2\pi^2}}
\sum_{k=0}^\infty {{(-1)^k} \over {k!}}\left[\gamma_k\left({1 \over 3}\right)-\gamma_k\left({2 \over 3}\right)\right],$$
$${4 \over {27}}\prod_{\overset{p \equiv 2}{mod ~3}} {{p^2-1} \over {p^2+1}}={1 \over 9}
\sum_{k=0}^\infty {{(-1)^k} \over {k!}}\left[\gamma_k\left({1 \over 3}\right)-\gamma_k\left({2 \over 3}\right)\right],$$
(b) (1.4b) holds, as well as
$$\prod_{\overset{p \equiv 1}{mod ~4}} {{p^2+1} \over {p^2-1}}={{12} \over \pi^2}G={4 \over {3\pi^2}} \sum_{k=0}^\infty {1 \over {k!}}\left[\gamma_k\left({1 \over 4}\right)-\gamma_k\left(
{3 \over 4}\right)\right],$$
(c)
$${1 \over \sqrt{5}}=\prod_{\overset{p \equiv 2}{mod ~5}} {{p^2-1} \over {p^2+1}}
\prod_{\overset{p \equiv 3}{mod ~5}} {{p^2-1} \over {p^2+1}}$$
$$=\sum_{k=0}^\infty {{(-1)^k} \over {k!}}\left[\gamma_k\left({1 \over 5}\right)-
\gamma_k\left({2 \over 5}\right)-\gamma_k\left({3 \over 5}\right)+
\gamma_k\left({4 \over 5}\right)\right],$$
$${{124} \over {125}}\zeta(3)= \prod_{\overset{p \equiv 2}{mod ~5}} {{p^3-1} \over {p^3+1}}
\prod_{\overset{p \equiv 3}{mod ~5}} {{p^3-1} \over {p^3+1}}$$
$$={1 \over {125}}\sum_{k=0}^\infty {{(-2)^k} \over {k!}}\left[\gamma_k\left({1 \over 5}\right)-
\gamma_k\left({2 \over 5}\right)-\gamma_k\left({3 \over 5}\right)+
\gamma_k\left({4 \over 5}\right)\right]$$
$$={1 \over {125}}\left[\psi''\left({2 \over 5}\right)-
\psi''\left({1 \over 5}\right)+\psi''\left({3 \over 5}\right)-
\psi''\left({4 \over 5}\right)\right],$$
(d) 
$${{3\pi^3} \over {64\sqrt{2}}}= {7 \over 8}\zeta(3)\prod_{\overset{p \equiv 5}{mod ~8}} {{p^3-1} \over {p^3+1}} \prod_{\overset{p \equiv 7}{mod ~8}} {{p^3-1} \over {p^3+1}}$$
$$={1 \over {8^3}}\sum_{k=0}^\infty {{(-2)^k} \over {k!}}\left[\gamma_k\left({1 \over 8}\right)+
\gamma_k\left({3 \over 8}\right)-\gamma_k\left({5 \over 8}\right)-
\gamma_k\left({7 \over 8}\right)\right],$$
$${{\pi^2} \over {\sqrt{2}}}= 6\zeta(2)\prod_{\overset{p \equiv 3}{mod ~8}} {{p^2-1} \over {p^2+1}} \prod_{\overset{p \equiv 5}{mod ~8}} {{p^2-1} \over {p^2+1}}$$
$$={1 \over {8}}\sum_{k=0}^\infty {{(-1)^k} \over {k!}}\left[\gamma_k\left({1 \over 8}\right)-
\gamma_k\left({3 \over 8}\right)-\gamma_k\left({5 \over 8}\right)+
\gamma_k\left({7 \over 8}\right)\right],$$
and (e)
$$\sum_{k=0}^\infty {{(-1)^k} \over {k!}}\left[\gamma_k\left({1 \over 8}\right)+
\gamma_k\left({3 \over 8}\right)-\gamma_k\left({5 \over 8}\right)-
\gamma_k\left({7 \over 8}\right)\right]
={{45} \over {32}}{{\zeta(4)} \over {\sqrt{2}G}}\prod_{\overset{p \equiv 7}{mod ~8}} \left({{p^2-1} \over {p^2+1}}\right)^2,$$
$$\sum_{k=0}^\infty {{(-2)^k} \over {k!}}\left[\gamma_k\left({1 \over 8}\right)-
\gamma_k\left({3 \over 8}\right)-\gamma_k\left({5 \over 8}\right)+
\gamma_k\left({7 \over 8}\right)\right]
={{3\pi^6} \over {1792\sqrt{2}}}{1 \over {\zeta(3)}}\prod_{\overset{p \equiv 7}{mod ~8}} \left({{p^3-1} \over {p^3+1}}\right)^2.$$


Let $E_j$ denote the Euler numbers, 
for which the initial values are $E_0=1$, $E_2=-1$, $E_4=5$, $E_6=-61$, $E_8=1385$,
$E_{10}=-50521$, and $E_{2n+1}=0$ for $n \geq 0$.

The next result complements Proposition 2 in the case that $D=-4$.
\newline{\bf Proposition 8}  
Let $k \geq 0$.  Then
$${{(-1)^k} \over {2(2k)!}}E_{2k}\left({\pi \over 2}\right)^{2k+1}=\sum_{n=0}^\infty {{(-1)^n} 
\over {n!}}\left[\gamma_n\left({1 \over 4}\right)-\gamma_n\left({3 \over 4}\right)\right](2k)^n.$$

{\bf Proposition 9}  
Sierpinski's constant \cite{sierpref}
$S=\gamma_{-4}\simeq 0.8228252496$ has the representation $S=2\gamma+{4 \over \pi}J_2$, wherein
$$J_2 \equiv \int_0^1 {{\ln(-\ln x)} \over {1+x^2}}dx$$
has the novel representation for $0<b<\pi/2$
$$J_2={1 \over 2}\sum_{k=1}^\infty {{E_{2k}b^{2k+1}} \over {(2k+1)!}}\left[\ln b -{b \over {(2k+1)}}\right]+{b \over 2}(\ln b-1)+{1 \over 2}\int_b^\infty {{\ln u} \over {\cosh u}}du.$$
In particular,
$$J_2=-{1 \over 2}\left[1+\sum_{k=1}^\infty {{E_{2k}} \over {(2k+1)!}}{1 \over {(2k+1)}}\right]
+{1 \over 2}\int_1^\infty {{\ln u} \over {\cosh u}}du.$$

{\bf Proposition 10}.  Let $\sigma_1(n)$ be the sum of divisors function and $K_0(z)$ the
zeroth order modified Bessel function of the second kind.  Then
(a)
$$\zeta\left({1 \over 2}\right)L_{-11}\left({1 \over 2}\right)=\gamma+\ln\left({\sqrt{11} \over
{8\pi}}\right)+4\sum_{n=1}^\infty (-1)^n \sigma_1(n)K_0(\sqrt{11}n\pi)$$
$$=\zeta\left({1 \over 2}\right){1 \over \sqrt{11}}\sum_{n=0}^\infty {1 \over 2^n}{1 \over {n!}}
\sum_{k=1}^{10}\left({k \over {11}}\right)\gamma_n\left({k \over {11}}\right),$$
and, more generally, (b) for $p >7$ a prime and class number $h(-p)=1$,
\footnote{so that $p>7$ is given by $11$, $19$, $43$, and $67$.}
$$\zeta\left({1 \over 2}\right)L_{-p}\left({1 \over 2}\right)=\gamma+\ln\left({\sqrt{p} \over
{8\pi}}\right)+4\sum_{n=1}^\infty (-1)^n \sigma_1(n)K_0(\sqrt{p}n\pi)$$
$$=\zeta\left({1 \over 2}\right){1 \over \sqrt{p}}\sum_{n=0}^\infty {1 \over 2^n}{1 \over {n!}}
\sum_{k=1}^{p-1}\left({k \over p}\right)\gamma_n\left({k \over p}\right).$$

{\bf Proposition 11}.  Let $0<s<1$.  Then $\zeta(s)<0$.  In particular, $\zeta(1/2)<0$.
In fact we also show that 
$$-{3 \over 2}+{1 \over {15 \sqrt{5}}}<\zeta\left({1 \over 2}\right)<-{{35} \over {24}}.$$

\medskip
\centerline{\bf Proof of Propositions}

{\it Proposition 1}.  (a) This follows from (1.6) and then (1.2),
$$L_{\pm k}(1)=\sum_{n-1}^\infty {{\chi_k(n)} \over n}={1 \over k} \sum_{m=1}^k \chi_k(m)
\zeta\left(1,{m \over k}\right)$$
$$={1 \over k}\sum_{m=1}^k\chi_k(m)\gamma_0\left({m \over k}\right),$$
(b) follows since $\sum_{n=1}^k \chi_k(n)=0$ for nonprincipal characters $\chi_k$.

{\it Proposition 2}. This again follows from (1.6) and the use of (1.2).

{\it Remark}.  From the relation $L_{-7}(2)=I_{-7}$ and \cite{coffeyjmp} we have for instance
$${1 \over \sqrt{7}}\sum_{n=0}^\infty{{(-1)^n} \over {n!}}\left[\gamma_n\left({1 \over 7}\right)+
\gamma_n\left({2 \over 7}\right)-\gamma_n\left({3 \over 7}\right)+
\gamma_n\left({4 \over 7}\right)-\gamma_n\left({5 \over 7}\right)-
\gamma_n\left({6 \over 7}\right)\right]$$
$$=4\left[3\mbox{Cl}_2(\theta_7)-3\mbox{Cl}_2(2\theta_7)+\mbox{Cl}_2(3\theta_7)\right].$$

{\it Proposition 3}.  The Dedekind zeta function of $\mathbb{Q}(\sqrt{-23})$ is given by
$$\zeta_{\mathbb{Q}(\sqrt{-23})}(s)=\zeta_0(s)+2\zeta_1(s)$$
$$=\sum_{m,n}'\left[{1 \over 2}{1 \over {(m^2+mn+6n^2)^s}}+{1 \over {(2m^2+mn+3n^2)^s}}\right],$$
with the prime on the sum indicating the exclusion of the term $m=n=0$.
We apply (1.2), (1.6), and the values \cite{zagiergangl}
$$\zeta_{\mathbb{Q}(\sqrt{-23})}(2)={{4\pi^2/3} \over {23^{3/2}}}[21D(z_1)+7D(z_2)+D(z_3)-3D(z_4)+D(z_5)]$$
$$=\zeta_0(2)+2\zeta_1(2),$$
where the partial zeta function values are given by
$$\zeta_0(2)=\zeta(4)+{{8\pi} \over {23^{3/2}}}+\left[{{\zeta(3)} \over 2}+\sum_{n=1}^\infty
\sigma_{-3}(n)(1+\pi n\sqrt{23})e^{-\pi n\sqrt{3}}\right],$$
and
$$\zeta_1(2)={{\zeta(4)} \over 4}+{{16\pi} \over {23^{3/2}}}+\left[{{\zeta(3)} \over 2}+\sum_{n=1}^\infty \sigma_{-3}(n)(1+\pi n\sqrt{23})e^{-\pi n\sqrt{3}}\right].$$

{\it Remark}.  Proposition 3 gives an example of the general relation for volumes in hyperbolic
$3$-space $\mathbb{H}^3$
$$\zeta_F(2)={{4\pi^2} \over D^{3/2}}\mbox{Vol}\left(\mathbb{H}^3/SL_2({\cal{O}}_F)\right),$$
where $F$ is an imaginary quadratic field.

{\it Proposition 4}. The results follow from (1.6) and the use of (1.2).

{\it Remark}.  In particular, we have corresponding summations of differences of Stieltjes
constants for values including
$$L_{-3}(3)={{4 \pi^3} \over {81\sqrt{3}}}, ~~~~L_{-4}(3)={\pi^3 \over {32}},$$
$$L_{-7}(3)={{32 \pi^3} \over {343\sqrt{7}}}, ~~~~L_{-8}(3)={{3\pi^3} \over {64\sqrt{2}}}.$$
For example,
$${{4 \pi^3} \over {81}}={1 \over 3^{5/2}}\sum_{n=0}^\infty{{(-1)^n} \over {n!}}2^n\left[\gamma_n\left({1 \over 3}\right)-\gamma_n\left({2 \over 3}\right)\right],$$
and
$${{32 \pi^3} \over {343}}={1 \over 7^{5/2}}\sum_{n=0}^\infty{{(-1)^n} \over {n!}}2^n\left[\gamma_n\left({1 \over 7}\right)+
\gamma_n\left({2 \over 7}\right)-\gamma_n\left({3 \over 7}\right)+
\gamma_n\left({4 \over 7}\right)-\gamma_n\left({5 \over 7}\right)-
\gamma_n\left({6 \over 7}\right)\right].$$

Therefore we have identities for the summation of differences of Stieltjes constants as
given by the following example.
\newline{\bf Corollary 1}.
$$\pi^3={{81} \over {4\cdot 3^{5/2}}}\sum_{n=0}^\infty{{(-1)^n} \over {n!}}2^n\left[\gamma_n\left({1 \over 3}\right)-\gamma_n\left({2 \over 3}\right)\right]$$
$$={{343} \over {32 \cdot 7^{5/2}}}\sum_{n=0}^\infty{{(-1)^n} \over {n!}}2^n\left[\gamma_n\left({1 \over 7}\right)+
\gamma_n\left({2 \over 7}\right)-\gamma_n\left({3 \over 7}\right)+
\gamma_n\left({4 \over 7}\right)-\gamma_n\left({5 \over 7}\right)-
\gamma_n\left({6 \over 7}\right)\right].$$

Another systematic way to produce summatory relations for differences of the Stieltjes
constants is to use the author's integral representation \cite{coffeyanly} 
for Re $a>0$,
$$\gamma_k(a)={1 \over {2a}}\ln^k a-{{\ln^{k+1} a} \over {k+1}}+{2 \over a}\mbox{Re}~
\int_0^\infty {{(y/a-i)\ln^k(a-iy)} \over {(1+y^2/a^2)(e^{2\pi y}-1)}}dy.$$
We then determine that, for instance,
$$\sum_{k=0}^\infty{{(-2)^k} \over {k!}}\left[\gamma_k\left({1 \over 3}\right)-\gamma_k\left({2 \over 3}\right)\right]={{243} \over {16}}$$
$$+54 \mbox{Re}\int_0^\infty\left[{1 \over {(3y+2i)^3}}-{1 \over {(3y+i)^3}}\right]{{dy} \over
{(e^{2\pi y}-1)}}dy.$$
For $q$ real, we observe that
$$\mbox{Re}~{1 \over {(3y+qi)^3}}=9y{{(3y^2-q^2)} \over {(9y^2+q^2)^3}}.$$
Then
$$\sum_{k=0}^\infty{{(-2)^n} \over {k!}}\left[\gamma_k\left({1 \over 3}\right)-\gamma_k\left({2 \over 3}\right)\right]={{243} \over {16}}
+456\int_0^\infty\left[{{3y^2-4} \over {(3y^2+4)^3}}-{{(3y^2-1)} \over {(3y^2+1)^3}}\right]{{ydy} \over {(e^{2\pi y}-1)}}.$$
By taking two derivatives of an integral representation for the digamma function, 
$$\psi(s)=\ln s-{1 \over {2s}}-2\int_0^\infty {{tdt} \over {(t^2+s^2)(e^{2\pi t}-1)}},$$
we obtain for the tetragamma function
$$\psi''(s)=-{1 \over s^2}-{1 \over s^3}+8\int_0^\infty {{(t^2-s^2)} \over {(t^2+s^2)^3}}{{tdt}
\over {(e^{2\pi t}-1)}}.$$
With the difference $\psi''(2/3)-\psi''(1/3)=8\pi^3/3^{3/2}$, we again obtain the result stated in the remark above.

{\it Proposition 5}. We use (1.2), (1.6), and the value
$$L_{-4}(4)=4^{-4}[\zeta(4,1/4)-\zeta(4,3/4)]={1 \over {1536}}\left[\psi'''\left({1 \over 4}
\right)-\psi'''\left({3 \over 4}\right)\right]$$
$$={1 \over 768}\psi'''\left({1 \over 4}\right)-{\pi^4 \over {96}}.$$
The latter relation follows from the functional equation 
$$\psi'''(1-z)-\psi'''(z)=-2\pi^4(2+\cos 2\pi z)\csc^4 \pi z.$$

{\it Proposition 6}. The Euler-Kronecker constant of $\mathbb{Q}(\sqrt{D})$ is defined as
$$\gamma_D=\gamma+{{L_D'(1)} \over {L_D(1)}}.$$
(a)  For $D=-4$ and $L_{-4}(1)=\pi/4$,
$$\gamma_{-4}=\gamma+{4 \over \pi}L_{-4}'(1)  
=\ln\left[2\pi e^{2 \gamma}{{\Gamma^2(3/4)} \over {\Gamma^2(1/4)}}\right],$$
Sierpinski's constant \cite{sierpref}.
With $L_{-4}(s)=4^{-s}[\zeta(s,1/4)-\zeta(s,3/4)]$,
$$L_{-4}'(s)=-(\ln 4)L_{-4}(s)+4^{-s}[\zeta'(s,1/4)-\zeta'(s,3/4)].$$
From (1.2),
$$\left[\zeta'(s,a)+{1 \over {(s-1)^2}}\right]_{s=1}=-\gamma_1(a),$$
giving
$$L_{-4}'(1)=-(\ln 4){\pi \over 4}+4^{-1}\left[\gamma_1\left({3 \over 4}\right)-
\gamma_1\left({1 \over 4}\right)\right],$$
and
$$\gamma_{-4}=\gamma-\ln 4+{1 \over \pi}\left[\gamma_1\left({3 \over 4}\right)-
\gamma_1\left({1 \over 4}\right)\right].$$
Now by the Kronecker limit formula \cite{ihara} we also have
$$\gamma_{-4}={\pi \over 3}-\ln 4+2 \gamma-4\sum_{k=1}^\infty \ln(1-e^{-2\pi k}).$$
By equating the two expressions for $\gamma_{-4}$ we obtain
$${1 \over \pi}\left[\gamma_1\left({3 \over 4}\right)-\gamma_1\left({1 \over 4}\right)\right]
={\pi \over 3}+\gamma-4\sum_{k=1}^\infty \ln(1-e^{-2\pi k}).$$
If we expand the logarithm, and interchange sums applying geometric series, we obtain
the stated result.

(b) For $D=-3$ we have
$$\gamma_{-3}=\gamma+{{L_{-3}'(1)} \over {L_{-3}(1)}}
=\gamma+{{3\sqrt{3}} \over \pi}L_{-3}'(1)  
=\ln\left[2\pi e^{2 \gamma}{{\Gamma^3(2/3)} \over {\Gamma^3(1/3)}}\right].$$
Since
$$L_{-3}'(s)=-(\ln 3)L_{-3}(s)+3^{-s}[\zeta'(s,1/3)-\zeta'(s,2/3)],$$
$$\gamma_{-3}=\gamma-\ln 3+{\sqrt{3} \over \pi}\left[\gamma_1\left({2 \over 3}\right)-
\gamma_1\left({1 \over 3}\right)\right].$$
Also
$$\gamma_{-3}={\pi \over {2\sqrt{3}}}-\ln 3+2\gamma -4\sum_{k=1}^\infty \ln |1-e^{-2\pi i\omega_{-3}k}|,$$
where $\omega_{-3}=(1+i\sqrt{3})/2$.
Again by expanding in the logarithmic sum and rearranging,
$$\sum_{k=1}^\infty \ln |1-e^{-2\pi i\omega_{-3}k}|=-\sum_{k=1}^\infty \sum_{\ell=1}^\infty
{{(-1)^{\ell k}} \over \ell} e^{-\sqrt{3}\pi k \ell}
=\sum_{\ell=1}^\infty {{(-1)^\ell} \over \ell}{1 \over {[(-1)^\ell-e^{\sqrt{3}\pi \ell}]}},$$
we obtain the stated result.

{\it Proposition 7}. 
(a) is demonstrated via (1.2) and the Euler product of $L_{-3}(s)=L(s,\chi_{-3})$, so that
$$\prod_{\overset{p \equiv 1}{mod ~3}} {{p^2+1} \over {p^2-1}}={{27L_{-3}(2)} \over {2\pi^2}},$$
and
$$\prod_{\overset{p \equiv 2}{mod ~3}} {{p^2+1} \over {p^2-1}}={{4\pi^2} \over {27L_{-3}(2)}}.$$

(b) is proved similarly using $L_{-4}(s)$ and (1.2), and in this case we exhibit more of the
details of the Euler products for this $L$ series and $\zeta(s)$.  The resulting prime products
are
$$\prod_{\overset{p \equiv 1}{mod ~4}} {{p^2+1} \over {p^2-1}}={{12G} \over \pi^2},$$
and
$$\prod_{\overset{p \equiv 3}{mod ~4}} {{p^2+1} \over {p^2-1}}={\pi^2 \over {8G}}.$$

We have
$$L_{-4}(s)=\prod_p (1-\chi_{-4}(p)p^{-s})^{-1}=\prod_{\overset{p \equiv 1}{mod 4}} (1-p^{-s})^{-1}
\prod_{\overset{p \equiv 3}{mod ~4}} (1+p^{-s})^{-1},$$
while
$$\zeta(s)=(1-2^{-s})^{-1}\prod_{\overset{p \equiv 1}{mod ~4}} (1-p^{-s})^{-1}
\prod_{\overset{p \equiv 3}{mod ~4}} (1-p^{-s})^{-1}.$$
Then
$$L_{-4}(s)=(1-2^{-s})^{-1}\zeta(s)\prod_{\overset{p \equiv 3}{mod ~4}} {{(1-p^{-s})} \over 
{(1+p^{-s})}}=(1-2^{-s})^{-1}\zeta(s)\prod_{\overset{p \equiv 3}{mod ~4}} {{(p^s-1)} \over 
{(p^s-1)}}.$$
Then $L_{-4}(2)=G$ and (1.2) are used.  

(c) uses (1.2), $L_5(s)=5^{-s}[\zeta(s,1/5)-\zeta(s,2/5)-\zeta(s,3/5)+\zeta(s,4/5)]$, 
its Euler product representation, 
$$L_5(s)=\prod_{\overset{p \equiv 1}{mod ~5}} (1-p^{-s})^{-1}
\prod_{\overset{p \equiv 2}{mod ~5}} (1+p^{-s})^{-1} \prod_{\overset{p \equiv 3}{mod 5}} (1+p^{-s})^{-1} \prod_{\overset{p \equiv 4}{mod ~5}} (1-p^{-s})^{-1},$$
and the relation of the latter to that of $\zeta(s)$, so that
$$L_5(s)={{\zeta(s)} \over {(1-5^{-s})}}\prod_{\overset{p \equiv 2}{mod ~5}} {{(1-p^{-s})} \over 
{(1+p^{-s})}}\prod_{\overset{p \equiv 3}{mod ~5}} {{(1-p^{-s})} \over {(1+p^{-s})}}.$$
Evaluations at $s=2$ and $s=3$ are then performed.

(d) uses (1.2), $L_{-8}(s)=8^{-s}[\zeta(s,1/8)+\zeta(s,3/8)-\zeta(s,5/8)-\zeta(s,7/8)]$, 
and $L_{+8}(s)=8^{-s}[\zeta(s,1/8)-\zeta(s,3/8)-\zeta(s,5/8)+\zeta(s,7/8)]$, 
their Euler product representations in relation to that of $\zeta(s)$, and respective evaluations 
at $s=3$ and $s=2$.

(e) evaluates $L_{-8}(2)$ and $L_{+8}(3)$ using (1.2) and the Euler products of $L_{\pm 8}(s)$
in terms of those of $L_{-4}(s)$, $\zeta(s)$, and $\zeta(2s)$.

{\it Proposition 8}.  This is an application of (1.2), (1.6), and the evaluation for $k\geq 0$
$$L_{-4}(2k+1)={{(-1)^k} \over {2(2k)!}}E_{2k}\left({\pi \over 2}\right)^{2k+1}.$$

{\it Proposition 9}.  Based upon the integral $J_2$ evaluated in the Discussion section of
\cite{coffeygammaratl}, we again find
$$S=\gamma-\ln 4+{1 \over \pi}\left[\gamma_1\left({3 \over 4}\right)-\gamma_1\left({1 \over 4}\right)\right].$$
Then we have the various integral representations
$$S=2\gamma+{4 \over \pi}\int_0^\infty {{e^{-x}\ln x} \over {1+e^{-2x}}}dx
=2\gamma+{2 \over \pi}\int_0^\infty {{\ln x ~dx} \over {\cosh x}}=2\gamma+{4 \over \pi}J_2.$$
By applying a generating function of the Euler numbers,
$$\mbox{sech} ~x=1+\sum_{k=1}^\infty {{E_{2k}} \over {(2k)!}}x^{2k}, ~~~~|x| < {\pi \over 2},$$
we obtain the stated form of $J_2$.

{\it Proposition 10}.  We outline the steps for the more general setting of part (b).
Under the hypotheses stated, according to \cite{selberg67} (p.\ 95)
$$\zeta\left({1 \over 2}\right)L_{-p}\left({1 \over 2}\right)=\gamma+\ln\left({\sqrt{p} \over
{8\pi}}\right)+2\sum_{n=1}^\infty (-1)^n \sigma_1(n)\int_0^\infty e^{-\pi n \sqrt{p}(y+1/y)/2} {{dy} \over y}$$
$$=\gamma+\ln\left({\sqrt{p} \over {8\pi}}\right)+4\sum_{n=1}^\infty (-1)^n \sigma_1(n)K_0(\sqrt{p}n \pi).$$
Then (1.2) and (1.6) are used.

{\it Remarks}.  Such values $L_p(1/2)$ are of particular interest since finding an instance
of $L_p(1/2)<0$ would be a violation of the generalized Riemann hypothesis. 

Of course from (1.2) we have
$$\zeta\left({1 \over 2}\right)=-2+\sum_{n=0}^\infty {\gamma_n \over {2^nn!}}$$
$$=-2+{1 \over \pi}\int_0^\infty [\ln(1+t)-\psi(1+t)]{{dt} \over t^{1/2}}
\simeq -1.46035 <0,$$
where $\gamma_n(1)=\gamma_n$.  Among many other occurrences of the value $\zeta(1/2)$, it
enters the Madelung constant $M_2=4(\sqrt{2}-1)\zeta(1/2)L_{-4}(1/2)$.  Hence we have
$$M_2=2(\sqrt{2}-1)\zeta\left({1 \over 2}\right)\sum_{n=0}^\infty {1 \over {2^n n!}}\left[
\gamma_n\left({1 \over 4}\right)-\gamma_n\left({3 \over 4}\right)\right].$$

{\it Proposition 11}.  For $0<\mbox{Re}~s<1$ \cite{debruijn}, 
$$\zeta(s)={1 \over {s-1}}+{{\sin\pi s} \over \pi}\int_0^\infty [\ln(1+t)-\psi(1+t)]{{dt} \over t^s}.$$
Now for $x>0$,
$${1 \over {2x}}<\ln x-\psi(x)<{1 \over x}.$$
Letting $0<s<1$,
$${1 \over 2}\int_0^\infty {{dt} \over {(1+t)t^s}}<\int_0^\infty [\ln(1+t)-\psi(1+t)]{{dt} \over t^s}<\int_0^\infty {{dt} \over {(1+t)t^s}}.$$
Then using a special case of the Beta function integral for Re$(a+b)>1$ and Re $b<1$,
$$\int_0^\infty {{dt}\over {(1+t)^a t^b}}=B(1-b,a+b-1)={{\Gamma(1-b)\Gamma(a+b-1)} \over
{\Gamma(a)}},$$
$${1 \over 2}<{{\sin\pi s} \over \pi}\int_0^\infty [\ln(1+t)-\psi(1+t)]{{dt} \over t^s}<1.$$
Then it follows that $\zeta(s)<0$.

We may considerably refine the last inequality by using the inequalities from \cite{gordon}
(Theorem 5)
$${1 \over {2x}}+{1 \over {12(x+1/4)^2}}<\ln x-\psi(x)<{1 \over {2x}}+{1 \over {12x^2}}.$$
Restricting for illustration to $s=1/2$,
$${1 \over 2}\int_0^\infty {1 \over t^s}\left[{1 \over {t+1}}+{1 \over {6(t+5/4)^2}}\right]dt
<\int_0^\infty [\ln(1+t)-\psi(1+t)]{{dt} \over t^{1/2}}<\int_0^\infty {1 \over t^{1/2}}
\left[{1 \over {t+1}}+{1 \over {6(t+1)^2}}\right]dt.$$
Thus we obtain
$${1 \over 2}\left(1+{2 \over {15\sqrt{5}}}\right)<{1 \over \pi}\int_0^\infty [\ln(1+t)-\psi(1+t)]{{dt} \over t^{1/2}} <{{13} \over {24}}.$$
Hence 
$$-{3 \over 2}+{1 \over {15 \sqrt{5}}}<\zeta\left({1 \over 2}\right)<-{{35} \over {24}}$$
follows.

{\it Remarks}.  Similarly we obtain 
$$-{5 \over 6}+{1 \over {\sqrt{2}\cdot 30\cdot 5^{1/4}}}<\zeta\left({1 \over 4}\right)=-{4 \over 3} +\sum_{n=0}^\infty \left({3 \over 4}\right)^n {\gamma_n \over {n!}} < -{{13} \over {16}},$$
which relates to $\zeta(3/4)$ via
$$\zeta\left({3 \over 4}\right)={{\sqrt{2+\sqrt{2}}\Gamma(1/4)} \over {(2\pi)^{1/4}}}\zeta\left(
{1 \over 4}\right).$$

From the value of $\zeta'(1/2)$ we obtain the identity
$${1 \over \pi}\int_0^\infty [\ln(1+t)-\psi(1+t)]{{\ln t} \over t^{1/2}}dt=\sum_{n=0}^\infty
{\gamma_{n+1} \over {2^n n!}}=-2\mbox{Re}~i\int_0^\infty {{\ln(1-iy)} \over {(1-iy)^{1/2}
(e^{2\pi y}-1)}}dy.$$ 





\pagebreak

\end{document}